\documentclass[english,letterpaper]{article}
\usepackage[english]{babel}
\usepackage{amsmath,amsthm,amssymb}

\usepackage{latexsym}
\usepackage{enumerate}
\usepackage{color}

\usepackage{amsmath}
\usepackage{amsfonts}
\usepackage{amsthm,amscd}
\usepackage{amsmath}
\usepackage{amssymb}
\usepackage{amstext}

\begin{document}

\title{ $C^{*}$- Algebras and  Thermodynamic Formalism}

\author{Ruy Exel (*) and Artur O. Lopes (**)}

\maketitle

%\begin{document}

\centerline{\bf S\~ao Paulo Journal of Mathematical Sciences - (USP-S\~ao Paulo)}

\centerline{\bf Vol. 2, 1 (2008), 285-307 - version updated in 2021}

%\vspace{0.5cm}
%\centerline{\bf}

\vspace{1.0cm}

\centerline{\bf *\, \,Departamento de Matem\'atica, UFSC,
Florian\'opolis, Brasil.}

\vspace{0.1cm}

\centerline{\bf and}

 \vspace{0.1cm}

\centerline{\bf **\,\,Instituto de Matem\'atica, UFRGS, Porto
Alegre, Brasil.}

\begin{abstract}

We  show a relation of the KMS state of a certain $C^{*}$-Algebra
${\cal U}$  with the Gibbs state of Thermodynamic Formalism. More
precisely, we consider here the shift $T:X \to X$ acting on the
Bernoulli space $X=\{1,2,...,k\}^{\bf N}$ and $\mu$ a Gibbs (equilibrium) state
defined by a Holder continuous normalized potential $p:X \to {\bf R}$, and
${\cal L}^2 (\mu)$ the associated Hilbert space.

Consider the  $C^{*}$-Algebra  ${\cal U}={\cal U}(\mu)$, which is
a sub-$C^{*}$-Algebra of the $C^{*}$-Algebra of linear operators
in   ${\cal L}^2 (\mu)$ which will be precisely defined later.
We call $\mu$ the reference measure.
Consider a fixed Holder potential $H>0$ and the  $C^{*}$-dynamical
system defined by the associated homomorphism $\sigma_t$. We are
interested in describe for such system the KMS states $\psi_\beta$
for all $\beta\in {\bf R}$.

We show a relation of a new Gibbs (eigenprobability of a Ruelle operator) probability $\nu_\beta$ to a KMS
state $\psi_{\nu_\beta}=\psi_\beta,$ in the $C^*$-Algebra ${\cal
U}={\cal U} (\mu)$, for every value $\beta\in {\bf R}$, where
$\beta$ is the parameter that defines the time evolution
associated to a homomorphism $\sigma_t=\sigma_{\beta i}$ defined
by the potential $H$. We show that for each real $\beta$ the KMS
state is unique and we explicit it.
The probability $\nu_\beta$ is the eigenprobability of the dual of the Ruelle operator of the non-normalized potential $
- \beta \log H$. The purpose of the present work is to explain (for an audience
which is more oriented to  Dynamical System Theory) part of the
content of a previous paper written by the authors.

\end{abstract}

\centerline {\bf Introduction} \vspace{1.0cm}

In this paper we  show a relation of the KMS state of a certain
$C^{*}$-Algebra ${\cal U}$ [BR] [P] [EL2] with the Gibbs state of
Thermodynamic Formalism [PP] [Bo] [R3]. The purpose of this work
is to explain for an audience which is more oriented to  Dynamical
System Theory part the content of the paper [EL3]. See also
[Re1],[Re2] for related material.

R. Bowen, D. Ruelle and Y. Sinai are the founders of what is
called in our days Thermodynamic Formalism Theory (see [PP] [R3]).

We will present initially the precise definitions we are going to consider.

We point out that we show here only the uniqueness part of the
results in [EL3]. The existence is based on the paper [W] which is
of Functional Analysis nature.

We refer the reader to [CL1] for a detailed analysis of the different meanings of the concept of Gibbs state from the point of view of Thermodynamic Formalism.

We consider here an expanding transformation $T:X \to X$ (to
simplify ideas one can consider the particular case where $T$ is
the shift acting on the Bernoulli space $X=\{1,2,...,k\}^{\bf
N}$). Consider $\mu$ the Gibbs state defined by a normalized Holder  continuous
potential $p:X \to {\bf R}$, and ${\cal L}^2 (\mu)$ the associated
Hilbert space. The function $p$ is sometimes called the Jacobian of $\mu$.

Consider the  $C^{*}$-Algebra  ${\cal U}={\cal U}(\mu)$, which is
a sub-$C^{*}$-Algebra of the $C^{*}$-Algebra of linear operators
in   ${\cal L}^2 (\mu)$ which will be precisely defined later.

We call $\mu$ the reference measure.

Consider a fixed Holder potential $H>0$ and the  $C^{*}$-dynamical
system defined by the associated $\sigma_z$. We are interested in
describe for such system the
KMS states $\psi_\beta$ for all $\beta\in {\bf R}$.

We show a relation of a new Gibbs probability $\nu_\beta$ to a KMS
state $\psi_{\nu_\beta}=\psi_\beta,$ in the $C^*$-Algebra ${\cal
U}={\cal U} (\mu)$, for every value $\beta\in {\bf R}$, where
$\beta$ is the parameter that defines the time evolution
associated to a homorphism $\sigma_t=\sigma_{\beta i}$ defined by
the potential $H$. We show that for each real $\beta$ the KMS
state is unique in Theorem 2.2. We present the explicit expression of $\psi_\beta$.

The probability $\nu_\beta$ is the Gibbs state (eigenprobability) for the potential
$ - \beta \log H$ (which is not normalized).

Given a potential $H$, we say the potential $\tilde{H}$ is cohomologous to
$H$, if there is $V$ such that $\log \tilde{H}= \log H - V + V\circ T $.

After we present our main results for Holder potentials in section
1 and 2  in section 3 we consider a non-Holder potential $H$ and we
will make an analysis  of phase transition nature (which do not
occur at $C^*$ Algebra level, in this case) associated to the KMS
problem in a case where $H$ can attain the value $1$ (and where
there is phase transition at Thermodynamic Formalism level).

Ground states in $C^*$-Algebras are also consider in the paper [EL3]. These corresponds to limits of the KMS state $\psi_\beta$ when $\beta \to \infty$.

An important contribution in the relation of $C^*$ algebras and Thermodynamic Formalism appears in chapters I.3 e II.5 in J. Renault Phd thesis  [Re] (se also [R0])

We refer the reader to [CL2] for a detailed description of phase transition in the sense of Thermodynamic Formalism.

\vspace{1.0cm}
\centerline {\bf Section 1 - KMS and Gibbs states}
\vspace{0.5cm}

We denote $C(X)$ the space of continuous functions on $X$
taking values on the complex numbers where $(X,d)$ is a compact metric space.

Consider the Borel sigma-algebra ${\cal B}$ over $X$ and
a continuous transformation $T:X \to X$. Denote by
${ \cal M} (T)$ the set of invariant probabilities for $T$.
We assume that $T$ is an expanding map.

We refer the reader to [Bo] [R1] [R2] [R3] [L4] for general definitions and properties
of Thermodynamic Formalism and expanding maps.

Tipical examples of such transformations (for which that are a lot
of nice results [R2]) are the shift in the Bernoully space and
also  $C^{1+\alpha}$-tranformations of the circle such that $|T
'(x)|> c > 1$, where $| \, \, \, |$ is the usual norm (one can
associate the circle to the interval $[0,1)$ in a standard way)
and $c$ is a constan.

The geodesic flow in compact constant negative curvature surfaces
induces in the boundary of Poincar\'e disk a Markov transformation
$G$ such that for some $n$, we have $G^n=T$, and where $T$ is
continuous expanding and acts on the circle (see [BS]). Our
results can be applied for such $T$.

We denote by ${\cal H}= {\cal H}_\alpha$ the set of
$\alpha$-Holder functions taking complex values, where
$\alpha$ is fixed $0< \alpha\leq 1$.

For each $\nu\in { \cal M} (T)$, the real non-negative value  $h(\nu)$ denotes the
the Shanon-Kolmogorov entropy of $\nu$ and
$h(T)= \sup \{h(\nu)|  \nu \in { \cal M} (T)\}$.
$h(T)$ is called the topological entropy
of $T$.

Given a continuous function $A: X \to \mathbb{R}$ we denote the Ruelle operator by ${\cal L}_A $ (which acts on continuous function $f$). More precisely if $g={\cal L}_A (f)$, then
$g(x)={\cal L}_A (f)(x) = \sum_{T(z)=x} e^{A} (z) f(z)$.

We say that the potential $A$ is normalized if ${\cal L}_A (1)=1$.

Given $A$, the dual operator ${\cal L}_A^* $ acts on probabilities on $\mathcal{M}(X)$.

We say that ${\cal L}_A^* (\nu)=\rho$ if for any continuous function $f$
$$\int f\, {\cal L}_A^* (\nu)\,=\,\int f \rho= \int {\cal L}_A (f)\, d \nu.$$

We denote by $\mu$ a fixed Gibbs state for a real
Holder potential $\log p: X\to
{\bf R}$. We suppose $\log p$ is already normalized [Bo][R3], in the sense that,
if ${\cal L}_{\log p}$ (for short
${\cal L}_p$) denotes the Ruelle-Perron-Frobenius operator for $\log p$,
that is for any $f:X\to {\bf C}$, and  all $x\in X$, we have
$({\cal L}_p (f))(x) = \sum_{T(z)=x} p(z) f(z)$, then we assume
that
${\cal L}_p (1)(x) = \sum_{T(z)=x} p(z) =1$ and $ {\cal L}_p^{*} (\mu)
= \mu.$

We will show later that the index $\lambda(x)=p(x)^{-1}$ for the
$C^{*}$-algebra associated to $\mu$.

As an interesting example we mention the case
where
$T$ has  degree
$k$, that is, for each $x\in X$ there exists
exactly $k$ different solutions $z$ for $T(z)= x$. We call each such $z$ a
pre-image of $x$.

If $T$ has degree $k$ and in the particular case where
$\mu$ is the maximal entropy measure (that is, $h(\mu)= h(T)=\log k$), then
$p=1/k$.

In order to simplify the arguments in our proofs we will assume
from now on that $T$ has
degree $k$.

One can consider alternatively in Thermodynamic Formalism ${\cal
L}_p$ acting on $C(X)$ or on ${\cal H}_\alpha$. Different spectral
properties for ${\cal L}_p$ ocurr in each one of these two cases
(see[Bo][R2]).

We will consider in the sequel a fixed real Holder-continuous
positive potential $H:X\to {\bf R}$ and ${\cal L}_{H,\beta}$,
$\beta \in {\bf R}$ the Ruelle-Perron-Frobenius operator for $ -
\beta \log H$, that is, for each continuous $f$ we have by
definition
$${\cal L}_{- \beta\, \log H} (f) (x)={\cal L}_{H,\beta} (f) (x)=
\sum_{T(z)=x} H(z)^{-\beta}  f(z)  .$$

We denote by $\lambda_{H,\beta}\in {\bf R}$
the largest eigenvalue of ${\cal L}_{H,\beta}$.
We also denote $\nu_{H,\beta}$
the unique probability  such that ${\cal L}_{H,\beta}^{*} (\nu_{H,\beta}) =
\lambda_{H,\beta}
\nu_{H,\beta}$,
and $h_{H,\beta}$ the unique function $h\in C(X)$ such that
$\int h d \nu_{H,\beta} =1 $ and ${\cal L}_{H,\beta} (h) = \lambda_{H, \beta}    h.$

As $H$ is fixed for good in order to simplify the notation  we will sometimes write
${\cal L}_{\beta}, {\cal L}_{\beta}^{*},
\lambda_{\beta},
\nu_{\beta}$, $h_{\beta}$.

$h_{\beta}$ is a real positive  Holder function.

The hypothesis about $H$ and $p$ being Holder in the Statistical
Mechanics setting means that in the Bernoulli space the
interactions between spins in neighborhoods positions decrease
very fast [L2] [L3]. In section 2.3 we will consider a non-Holder
potential $H$ where in this case it will appear a phase-transition
phenomena. This model is known as the Fisher-Felderhof model [FF],
[L2], [L3], [FL]. In this case the interactions do not decrease so
fast.

We return now to the Holder case.

It is well known the variational principle for such potential
$ - \beta \log H$,
$$ P_H (\beta) = \log \lambda_{H,\beta}= \sup \{  h(\nu) +
\int (-
\beta \log H) d \nu | \nu \in {\cal M}(T)\}  .$$

The probability $\mu_{H,\beta}=h_{H,\beta} \nu_{H,\beta} \in {\cal M} (T)$ and satisfies
$$\sup \{  h(\nu) +
\int ( - \beta \log H) d \nu| \nu \in {\cal M}(T)\}
=$$
$$ h(\mu_{H,\beta}) +
\int ( - \beta \log H) d \mu_{H,\beta}.$$

{\bf Definition 1.1:}
The probability $\mu_\beta =h_{H,\beta} \nu_{H,\beta} $ is called equilibrium state
for the function $ - \beta \log H$
where $\beta$ and $H$ are fixed.

{\bf Definition 1.2:}
The probability $\nu_{H,\beta}$ is called eigenmeasure or  Gibbs state for the
function $  - \beta \log H$
where $\beta$ and $H$ are fixed.  It satisfies
$${\cal L}_{H,\beta}^{*} (\nu_{H,\beta})= \lambda_{H,\beta}\nu_{H,\beta}.$$

The probability
$ \mu_{H,\beta}$ is unique for the variational problem and
$\nu_{H,\beta}$ is
unique for the the eigenmeasure problem associated to
the value $\lambda_{H,\beta}$,
if $p$ and $H$ are Holder. If we do not assume $p$ and $H$ Holder
then there exist counterexamples for uniqueness in both cases [L2] [L3].
We will return to this point later.

For some reason the eigen-probabilities have a distinguished role
here, but not the equilibrium states.

$P_H (\beta)$ is called the pressure of $ - \beta \log H$ (or
sometimes Free-Energy) and is a convex analytic function of
$\beta$.

If $T$ has degree $k$ and in the particular case where
$\mu$ is the maximal entropy measure (that is, $h(\mu)= h(T)=\log k$), then
$p=1/k$.

\medskip

We consider the $C^*$-Algebra $
L({\cal L}^2 (\mu))$
of bounded linear operators acting on  ${\cal L}^2 (\mu) $ with the strong norm. The operation $*$ on operators is the one induced from the inner product on  ${\cal L}^2 (\mu)$.

\medskip

{\bf Definition 1.3:}
Denote by $S: {\cal L}^2 (\mu) \to {\cal L}^2 (\mu)$ the Koopman  operator
where for $ \eta \in {\cal L}^2 (\mu)$ we define
$(S \eta) (x) = \eta (T(x))$. Such $S$ defines a linear bounded operator in
${\cal L}^2 (\mu)$.

In Thermodynamic Formalism
it is usual to consider the Koopman operator  acting on ${\cal L}^2 (\mu)$
(the space of complex
square integrable functions over ${\cal L}^2 (\mu)$), and it is well known that
its adjoint (over ${\cal L}^2 (\mu)$) is the operator
${\cal L}_p=S^*$ acting on ${\cal L}^2 (\mu)$.

As we assume $X$ is compact, any continuous function $f$ is in ${\cal L}^2 (\mu)$.

{\bf Definition 1.4:}
Another important class of linear operators is
$M_f : {\cal L}^2 (\mu) \to {\cal L}^2 (\mu)$, for a given fixed $f\in C(X)$,
and defined by $M_f (\eta) (x)= f(x) \eta (x)$, for any $\eta$ in ${\cal L}^2 (\mu)$.

In order to simplify the notation, sometimes
we denote by $f$ the linear operator $M_f$.

Note that for $M_f $ and $M_g$, $f,g \in C(X)$,  the product operation satisfies
$M_f \circ M_g = M_ {f.g}$, where $\,.\,$ means multiplication
over the complex field ${\bf C}$.

Note that the   $*$ operation applied on $M_f$, $f \in C(X)$,
is given by $M_f^{*}= M_{\overline{f}}$, where $\overline{z}$
is the complex conjugated
of  $z\in{\bf C}$.
In this sense, $M_f^{*}$ is the adjoint operator of $M_f$ over
${\cal L}^2 (\mu)$.

The main point for our choice of $\mu$ as eigen-probability for
${\cal L}_p^*$, is that in ${\cal L}^2 (\mu)$,  the dual of the
Koopman operator $S$ is the operator ${\cal L}_p=S^*$ acting on
${\cal L}^2 (\mu)$. Indeed, for any $f,g$ we have
$$ \int\, f \, (g\circ T)\, d \mu\, =\,   \int\, f \, (g\circ T)\, d{\cal L}_p^*( \mu)=
\int\,{\cal L}_p(\, f \, (g\circ T)\,)\, d \mu= \int\,{\cal
L}_p(\, f \,)\,  g\, d \mu.$$

\medskip

It is important not confuse the dual of the Ruelle operator ${\cal L}_p$ in the Hilbert structure sense with the dual of
${\cal L}_p$ as a linear functional on continuous functions.

\medskip

$L ({\cal L}^2 (\mu))$, the set of linear operators over ${\cal
L}^2 (\mu)$, is a very important $C^{*}$-Algebra. We will analyze
here a sub-$C^{*}$-Algebra of such $C^{*}$-Algebra (defined with
the above operations  $.$ and $*$), more precisely the
$C^{*}$-Algebra ${\cal U}$.

{\bf Definition 1.5:} We denote by $\alpha: C(X) \to C(X)$  the
linear operator such that for any $f$, we have $\alpha (f)= f\circ
T$.

We have to show how the operators $S$ and $M_f$ acting on ${\cal
L}^2 (\mu)$ interact with the operators ${\cal L}_p$ and $\alpha$
acting on $C(X)$.

One can easily see that $\alpha(M_f)=M_{f\circ T}$. This is
the first relation.

In  the simplified  notation (we identify
$M_f $  with $f$), one can read last expression as
$\alpha(f)= f \circ T$.

In this way $\alpha^n(f)= f \circ T^n$.

\medskip

If $\cal{B}$ is the Borel sigma-algebra then we denote by $\mathcal{F}_n$ the Sigma-algebra $T^{-n} \cal{B}$.

It is know that if we consider the  probability  $\mu$, then the conditional expected value
$$ E( f \,| \,\mathcal{F}_n)=  E_\mu( f \,| \,\mathcal{F}_n)=\alpha^n ( \mathcal{L}^n_p (f)). $$

More precisely
\begin{equation} \label{ah}  E( f \,| \,\mathcal{F}_n) (x) = \mathcal{L}^n (f) \, (\sigma^n (x)).
\end{equation}

As $\mathcal{F}_m \subset \mathcal{F}_n$ for $ m\geq n$,  we have
$$ \,  E_\mu (\,(\,  E_\mu\,(f | \,\mathcal{F}_m\,)\,)\, | \,\mathcal{F}_n\,) = E_\mu\,(f | \,\mathcal{F}_m\,),$$
and
$$ \,  E_\mu (\,(\,  E_\mu\,(f | \,\mathcal{F}_n\,)\,)\, | \,\mathcal{F}_m\,) = E_\mu\,(f | \,\mathcal{F}_m\,).$$

\medskip

{\bf Definition 1.6:} Consider  the $C^{*}$-Algebra contained in the set of bounded operators $
L({\cal L}^2 (\mu))$ generated by the elements of the form
$M_f S^n (S^{*})^n M_g$, where $n\in {\bf N}$ and $f,g\in C(X)$.
We denote such $C^{*}$-Algebra by ${\cal U}= {\cal U} (\mu,T).$ We
call ${\cal U}$ the  $C^{*}$-Algebra associated to $\mu$.

Each element $a$ in ${\cal U}$ is the limit of finite
sums  $\sum_{i} M_{f_i} S^{n_i} (S^{*})^{n_i} M_{g_i}$.

$C(X)$ is contained in $\mathcal{U}$, via $M_f$, where $f$ is any continuous function $f:X \to \mathbb{R}.$

Note that $f\to M_f$ defines a linear injective function of $C(X)$
on ${\cal U}$.

\medskip
We denote $e_n = S^n \, (S^n)^*= E( f \,| \,\mathcal{F}_n) \in \mathcal{U}.$

\bigskip

{\bf Important Properties:}

\medskip

We have  basic relations in such $C^{*}$-Algebra ${\cal U}$:
\medskip

a) $ (S^{*})^n S^n  =1 $, for all $n\in {\bf N}$ (it follows from
$ S^{*} S =1 $) .

proof: for any $\eta \in {\cal L}^2(\mu)$, we have
$$ S^* \, S (\eta)(x) = {\cal L}_p (\eta (T(.))(x) =\sum_{T(y)=x}
\, p(y) \, \eta(T(y)) =\sum_{T(y)=x} \, p(y) \, \eta(x)\,=\,
\eta(x).$$

That is, $ (S^{*})^n S^n  $ is the identity operator.

b) $(S^{*})^n M_f S^{n} = M_{{\cal L}_p^n (f)}$, for all $n \in
{\bf N}$, $f\in C(X)$ (it follows from $S^{*} M_f S= M_{{\cal L}_p
(f)}$).

proof: for any $\eta \in {\cal L}^2(\mu)$, we have
$$ S^{*} \, M_f \,S \, (\eta) \, (x)   = {\cal L}_p (f\, \eta(T(.)))(x)=
{\cal L}_p (f) \, (x)\, \eta (x) .   $$

c) $S M_f = \alpha (f) S$ for any continuous $f$, that is, for any
$\eta \in {\cal L}^2 (\mu)$, $ S M_f (\eta) = f\circ T . \eta
\circ T= \alpha (f) . S  (\eta).$

d) $[e^n \, M_f ](\eta)= [\,S^n\, (S^*)^n M_f\,]  (\eta) = E_\mu (\,(f\, \eta)\,| \,\mathcal{F}_n\,)$.

e) $  e^n (\eta)=\,S^n\, (S^*)^n   (\eta) = E_\mu (\,\, \eta\,| \,\mathcal{F}_n\,)$.

f) $M_f\, e^n (\eta)= [\,M_f \, S^n\, (S^*)^n ] (\eta) = f\,\, E_\mu (\,\eta\,| \mathcal{F}_n\,)$.

g) $M_f\, e^n \, M_g (\eta)=$
\begin{equation} \label{isto2} [\,M_f \, S^n\, (S^*)^n\, M_g ] (\eta) = f\,\, E_\mu (g\,\eta\,| \mathcal{F}_n\,)= f \,[ \,\mathcal{L}_p^n (g\, \eta\, ) \,(\sigma^n)\,].
\end{equation}

h) $ S^n\, (S^*)^n\, M_g   S^n\, (S^*)^n\,  =         \, E_\mu (g\,\,| \mathcal{F}_n\,) \, (\,S^n\, (S^*)^n\,)  \, E_\mu (g\,\,| \mathcal{F}_n\,)\, e^n    $ because
$$ [S^n\, (S^*)^n\, M_g   S^n\, (S^*)^n ](\eta)=  E_\mu (g\,\,   E_\mu (\eta\,\,| \mathcal{F}_n\,)      \, | \, \mathcal{F}_n\,) =     E_\mu (g\,\,   | \mathcal{F}_n\,)    \,\,  E_\mu (\eta\,\,| \mathcal{F}_n\,) .$$

i) If  $n\leq m$ we have
$$ [M_f\, e^n\, M_g   e^m\, M_h ]=[M_f\, S^n\, (S^*)^n\, M_g   S^m\, (S^*)^m\, M_h ]=$$
  $$[M_f\,(\, S^n\, (S^*)^n\, M_g   S^m\, (S^*)^m\,)\, M_h ]=      M_f     \, E_\mu (g\,\,| \mathcal{F}_n\,) \, (\,S^m\, (S^*)^m\,) \, M_h     =    $$
$$M_f     \, E_\mu (g\,\,| \mathcal{F}_n\,) \, e^m  \, M_h     .   $$

j) If  $n\geq m$ we have
$$[M_f\, e^n\, M_g   e^m\, M_h ]= [M_f\, S^n\, (S^*)^n\, M_g   S^m\, (S^*)^m\, M_h ]=  M_f \,e^n    \, E_\mu (g\,\,| \mathcal{F}_m\,) \,   \, M_h  $$

proof: note first that taking adjoint with respect to the $\mathcal{L}^2 (\mu)$ structure
$$ ( M_f\, S^n\, (S^*)^n\, M_g \,)^* = ( M_g\, S^n\, (S^*)^n\, M_f \,)   .  $$

Then,
$$ [M_f\, S^n\, (S^*)^n\, M_g   S^m\, (S^*)^m\, M_h ]^* = M_h\, S^m\, (S^*)^m\, M_g   S^n\, (S^*)^n\, M_f $$ and we can apply item i) to get
$$ [M_f\, S^n\, (S^*)^n\, M_g   S^m\, (S^*)^m\, M_h ]^* = M_h     \, E_\mu (g\,\,| \mathcal{F}_m\,) \, e^n  \, M_f. $$

Now taking adjoint once more we get
$$ [M_f\, S^n\, (S^*)^n\, M_g   S^m\, (S^*)^m\, M_h ]=  M_f \,e^n    \, E_\mu (g\,\,| \mathcal{F}_m\,) \,   \, M_h.$$

\vspace{0.3cm}

{\bf Example 1:}
$$ [\,M_f \, e^3 M_g\,   e^4 M_h \,] (\eta)\,(x)= $$
$$ [\,M_f \, S^3\, (S^*)^3 M_g\,   S^4\, (S^*)^4 M_h \,] (\eta)\,(x)= $$
$$  \,M_f \, S^3\, (S^*)^3 M_g\,[\,  E_\mu (\,(f\, \eta)\,| \,\mathcal{F}_4\,)(x)]=$$
$$  \,M_f \, S^3\, (S^*)^3 \,[\,g(x)\,  E_\mu (\,(f\, \eta)\,| \,\mathcal{F}_4\,)(x)]=$$
$$  \,f(x)  \,  E_\mu (\,\,\,g(x)\,\,\,\,  E_\mu (\,(f\, \eta)\,| \,\mathcal{F}_4\,)(x)\,\,\,|\,\mathcal{F}_3 \,)\,(x).$$

If $u$ is $\mathcal{F}_m $ measurable and $m>n$, then $u$ is $\mathcal{F}_n $ measurable.

Then,

$$ [\,M_f \, S^3\, (S^*)^3 M_g\,   S^4\, (S^*)^4 M_h \,] (\eta)\,(x)= $$
$$  \,f(x)  \, E_\mu (\,(f\, \eta)\,| \,\mathcal{F}_4\,)(x)\,\,\,  E_\mu (\,g(x)\, \,|\,\mathcal{F}_3 \,)\,(x).$$

By the other hand

$$ [\,M_f \, e^4 M_g\,   e^3 M_h \,] (\eta)\,(x)= $$
$$ [\,M_f \, S^4\, (S^*)^4 M_g\,   S^3\, (S^*)^3 M_h \,] (\eta)\,(x)= $$
$$  \,f(x)  \,  E_\mu (\,\,\,g(x)\,\,\,\,  E_\mu (\,(f\, \eta)\,| \,\mathcal{F}_3\,)(x)\,\,\,|\,\mathcal{F}_4 \,)\,(x).$$

 \vspace{0.2cm}

 {\bf Remark 0:} If we consider the $C^*$-algebra
generated $M_f \,S^m (S^{*})^n\, M_g$, where $n,m\in {\bf N}$ and
$f,g\in C(X)$, we have a different setting (which is  usually
called a  Vershik  $C^*$-algebra) which was consider in another
paper by R. Exel [E3]. In this case, the KMS state exists only for
one value of $\beta$.

\vspace{0.2cm}

We now return to our setting.

An extremely important result will be shown in expression (*1) and (*2) in Lemma 2.1 which
claims that there exists functions $u_i$, $i\in \{1,2,..,k\}$,
such that
$$\sum_{i=1}^k M_{u_i} S S^{*} M_{u_i} = 1.$$
\vspace{0.3cm}

A bijective linear transformation $ K: {\cal U} \to {\cal U}$ which preserves the
composition and the $*$ operation is called an automorphism  of ${\cal U}$

We denote by Aut(${\cal U}$) the set of automorphism of the
$C^{*}$-Algebra ${\cal U}$.

{\bf Definition 1.7:}
Given a positive function $H$ we define the group homomorphism
$\sigma_t$, where for each  $t\in {\bf R}$ we have
$\sigma_t \in $ Aut(${\cal U}$) [BP] [P], is defined by:

a) for each fixed $t\in  {\bf R}$ and any $M_f$, we have
$\sigma_t (M_f)=M_f$,

b) for each fixed $t\in  {\bf R}$, we have
$\sigma_t (S)= M_{H^{i  t}} \circ S$,
, in the sense that
$(\sigma_t  (S) (\eta)) (x) = H^{i t}(x) \eta(T(x))\in
{\cal L}^2 (\mu),$ for any $\eta\in
{\cal L}^2 (\mu)$.

The value $t$ above is related to temperature and not time, more
precisely we are going to consider bellow $t= \beta i$ where
$\beta $ is related to the inverse of temperature in Thermodynamic
Formalism (or Statistical Mechanics).

It can be shown that
for each $t$ fixed, we just have to define $\sigma_t$ over
the generators of ${\cal U}$ in order to define
$\sigma_t$ uniquely on ${\cal U}$. In this way a) and b)
above define $\sigma_t$.

We will assume in this section from now on that $H$ is Holder
in order we can use the strong results of Thermodynamic Formalism.

{\bf Remark 1:} Note that for $\eta\in {\cal L}^2 (\mu)$, we have
$$(\sigma_{t } (S^2) \eta) (x)=
\sigma_{t }  (M_{H^{i t }} (\eta \circ T) )(x)=
$$
$$
M_{H^{i t }}M_{H^{i t }\circ T} (\eta \circ T^2) (x),$$
therefore
$\sigma_{t } (S^2)= H^{t i } (H\circ T)^{t i } S^2.$
It follows easily by induction that
$$\sigma_{t } (S^n)=
\Pi_{j=0}^{n-1} (H\circ T^j )^{t i } S^n.$$

Taking dual in both sides of the above expression we get other important relation
$$\sigma_{t} ((S^{*})^n)=
(S^{*})^n \Pi_{j=0}^{n-1} (H \circ T^j )^{-t i }.$$

Finally,

\begin{equation*}
\sigma_{t}  (M_{f_2}
S^m (S^{*})^m M_{g_2}) = M_{f_2} H^{t i \,[m] }
S^m (S^{*})^m H^{-t \, i\,[m]}\,M_{g_2},  \,\,\,\,\,\,\,\,(*5)
\end{equation*}

where
 $H^{t\,i  \, [m]}(x)= \Pi_{i=0}^{m-1} H(T^i (x))^{t\,i }.$

From f) above we get for $t=i$
$$
\sigma_{i}  (M_{f_2}
S^m (S^{*})^m M_{g_2}) (\eta)(x) = $$
$$[\,M_{f_2} H^{i i \,[m] }
S^m (S^{*})^m H^{-\,i \, i [m]}\,M_{g_2}]\,(\eta)(x)=
$$
$$  [\,M_{f_2} H^{-\,[m] }
S^m (S^{*})^m H^{[m]}\, M_{g_2}\,]\,(\eta)(x)=
$$
$$ f_2(x)\, H^{\,-[m] }(x)\, E_\mu (H^{[m]}\,g_2\,\eta\,| \mathcal{F}_m\,)\,(x).$$

In terms of the formalism of  $C^{*}$-dynamical systems,
the positive function $H$ defines the dynamics
of the evolution with time $t\in {\bf R}$ of a  $C^{*}$-dynamical system.
Our purpose is to analyze such system  for each pair $(H,\beta)$.

{\bf Definition 1.8:} An element $a$ in a $C^{*}$-Algebra is positive, if it is of the
form $a=b \,b^{*}$ with $b$ in the $C^{*}$-Algebra.

{\bf Definition 1.9:} By definition a  "$C^{*}$-dynamical system state"
is a linear functional $\psi: {\cal U} \to {\bf C}$
such that

a) $\psi (M_1)=1$

b) $\psi (a) $ is a positive real number
for each positive element      $a$ on the $C^{*}$-Algebra
${\cal U}$.

A  "$C^{*}$-dynamical system state" $\psi$ in $C^{*}$-dynamical systems
plays the role of a
probability
$\nu$ in Thermodynamic Formalism. For a fixed
$H$, we have a dynamic temporal
evolution defined by $\sigma_t$ where $ t \in {\bf R}$.

{\bf Definition 1.10:}
An element $a\in {\cal U}$ is called analytic for $\sigma$ if
$\sigma_t(a)$ has an analytic extension from $t\in{\bf R}$ to all
$t\in {\bf C}$.

{\bf Definition 1.11:} For a fixed $\beta\in {\bf R}$ and $H$, by
definition, $\psi$ is a {\bf KMS state associated to $H$ and
$\beta$} in the  $C^{*}$-Algebra ${\cal U} (\mu,T)$, if  $\psi$ is
a $C^{*}$-dynamical system state, such that for any $b\in {\cal
U}$ and any analytic $a\in {\cal U}$ we have

$$ \psi (a . b) = \psi (b . \sigma_{\beta i} (a)).$$

For $H$ and $\beta$ fixed,
we denote a KMS state by $\psi_{H,\beta}=\psi_\beta$
and we leave $\psi$ for a general $C^{*}$-dynamical system state.

It is easy to see  that for $H$ and $\beta$ fixed, the condition
$$ \psi (a . b) = \psi (b . \sigma_{\beta i} (a)),$$
is equivalent to $\forall \tau\in {\bf C}$,
$$ \psi (\sigma_{\tau}(a) . b) = \psi (b . \sigma_{\tau + \beta i} (a)).$$

It follows from section 8.12 in [P] that if $\psi_\beta$ is a KMS state
for $H,\beta$, then for any analytic $a\in {\cal U}$, we
have that
$\tau\to \psi_\beta(\sigma_\tau ( a) ) $ is
a bounded entire function
and therefore constant. In this sense $\psi$ is stationary for the continuous time evolution defined by the flow $\sigma_t$.

Note that the KMS state, in principle, could depend of the initially
chossen $\mu$ because we  are considering  ${\cal L}^2 (\mu)$ when defining
${\cal U}$, but in the end it will be defined by a measure
that depends only in $\beta$ and $H$

We point out that it can be shown that in order to characterize
$\psi$ as a KMS  state we just have to check the condition
$ \psi (a . b) = \psi (b . \sigma_{\beta i} (a))$ for $a,b$ the linear generators
of ${\cal U}$, that is,
$a$ of the form
$M_{f_1} S^{n} (S^{*})^{n} M_{g_1}$
and
$b$ of the form
$M_{f_2} S^{m} (S^{*})^{m} M_{g_2}$.

A natural question is: for a given $\beta$ and $H$,
when the KMS state $\psi_{H,\beta}$ exist and when it  is unique?

We are interested mainly in uniqueness and explicitly. We will explain this point
more carefully later.

{\bf Remark 2:} Note that when $\psi$ is a KMS state,
$\psi (f\,. a\,. g)= \psi (\sigma_{\beta i} (g)\, f\, a)= \psi (g . f\,  a)= \psi (f . g\,  a)$, for
any    $f,g\in C(X)$ and $a\in {\cal U}$.

Our purpose here is to show how to associate in a unique way
each KMS state $\psi_{H,\beta}=\psi_\beta$ to the eigenmeasure $\nu_{H,\beta}= \nu_\beta$ defined before.

Remember that over  ${\cal L}^2 (\mu)$ the operator ${\cal L}_p=S^*$ is
adjoint
of the operator $f \to S(f)=f\circ T$.

We call $\lambda(x)= p(x)^{-1}$ the index and
we denote by
$$\lambda^{[n]} (x)= (p(x) p (T(x))...p(T^{n-1} (x)))^{-1}.$$

\medskip

We denote  $H^{\beta  \, [n]}(x)= \Pi_{i=0}^{n-1} H(T^i (x))^{\beta } $ and $\Lambda_n = H^{-\beta  \, [n]}\,\lambda^{[n]}.$

\medskip

From this follows that for any continuous function $f$ we have ${\cal L}^n_{\beta} (f ) = {\cal L}^n_{p} (\Lambda_n \,f  ).$

Remember that for any continuous function $k$ we have $\mathcal{L}^n_p (k \circ T^n)\,(x)= k(x)$ because $\mathcal{L}^n_p(1)=1$.

\medskip

{\bf Lemma 1.1}  For any any $\beta$ and continuous function $f$

\begin{equation} \label{op1}\int f\, d \nu_{\beta} = \int (\Lambda_n)^{-1} E_\mu( \Lambda_n\, f \,| \,\mathcal{F}_n\,) \,d \nu_{\beta}.
\end{equation}

{\bf Proof:}

Note that
$$\int {\cal L}^n_{p} (\Lambda_n \,f  )  \,d \nu_{\beta} =\int   {\cal L}^n_{\beta} ( \,f  ) d \nu_\beta =\lambda_\beta^n  \int f d \nu_\beta.$$

Now taking $f = (\Lambda_n)^{-1} \alpha^n (g)= (\Lambda_n)^{-1} \,( g \circ T^n)$ we get  from above

$$\int g\, d\nu_\beta= \int {\cal L}^n_{p} (\Lambda_n \, (\Lambda_n)^{-1} \,( g \circ T^n)  )  \,d \nu_{\beta}  =\lambda_\beta^n  \int  (\Lambda_n)^{-1} \,( g \circ T^n) d \nu_\beta.$$

Now taking $g= \mathcal{L}^n_p (\Lambda_n\,f)$ we get

$$ \int (\Lambda_n)^{-1} E_\mu( \Lambda_n\, f \,| \,\mathcal{F}_n\,) \,d \nu_{\beta} = \int (\Lambda_n)^{-1}  \alpha^n ( \mathcal{L}^n_p (\Lambda_n\,f))\,d \nu_{\beta} =$$
$$\lambda_\beta^{-n} \,\int\, \mathcal{L}^n_p (\Lambda_n\,f)\,d \nu_{\beta} = \lambda_\beta^{-n} \,\int\, \mathcal{L}^n_\beta (\,f)\,d \nu_{\beta}=\, \int f \, d \nu_\beta.$$

\qed

\bigskip

\vspace{1.0cm}
\centerline {\bf Section 2 - The main result}
\vspace{0.5cm}

We define $G: \mathcal{U} \to C(X)$ by $G(M_f e_n M_g)= f \lambda^{-[n]} \,g$ where $e_n = S^n (S^{*})^n.$

Moreover, $G (M_f  M_g)=f\, g$

Note that  we define $G$ in the elements of the form $M_f e_n M_g$, $n\geq 0,$ and then we define $G$ in  $\mathcal{U} $ by linear combinations and limits.

Suppose $\phi=\phi_{\nu}: C(X) \to \mathbb{C}$ is of the form $\phi(f) = \int f d \nu$ where $\nu$ is a probability on $X$.

There is a canonical way to define a $C^*$-dynamical  system  state $\psi_\nu : \mathcal{U} \to \mathbb{C} $ by
$$ \psi_\nu(M_f e_n M_g) = \phi_\nu (G(M_f e_n M_g))\,=\, \int f \lambda^{-[n]} \,g\, d \nu .  $$
\medskip

In this way if $n\leq m$ (by item i)\,)
$$  \psi_\nu (\,M_f\, S^n\, (S^*)^n\, M_g   S^m\, (S^*)^m\, M_h \,)=$$
  $$\psi_\nu (   M_f     \, E_\mu (g\,\,| \mathcal{F}_n\,) \, (\,S^m\, (S^*)^m\,) \, M_h   )=
  $$
  $$
\psi_\nu (   M_f     \, E_\mu (g\,\,| \mathcal{F}_n\,) \, e^m \, M_h   )= \int  \, E_\mu (g\,\,| \mathcal{F}_n\,)\,f\,h\, \lambda^{-[m]} \,\, d \nu . $$

In this way if $n\geq m$ (by item j)\,)
$$  \psi_\nu (\,M_f\, S^n\, (S^*)^n\, M_g   S^m\, (S^*)^m\, M_h \,)=$$
$$  \psi_\nu (
 M_f \,e^n    \, E_\mu (g\,\,| \mathcal{F}_m\,) \,   \, M_h \,)=\int  \, E_\mu (g\,\,| \mathcal{F}_m\,)\,f\,h\, \lambda^{-[n]} \,\, d \nu . $$

{\bf Theorem 2.1:} Given $\phi_\nu$ and $\psi_\nu = \phi_\nu \circ G$ we get that
$\psi_\nu$ is KMS for temperature $\beta$, if and only if, $\phi_\nu$ satisfies

$$\phi_\nu  (f\,) = \phi_\nu (\,\,(\Lambda_n)^{-1} E_\mu( \Lambda_n\, f \,| \,\mathcal{F}_n\,) \,), $$

which is the same that to say that $\nu$ satisfies

$$\int f\, d \nu = \int (\Lambda_n)^{-1} E_\mu( \Lambda_n\, f \,| \,\mathcal{F}_n\,) \,d \nu. $$

{\bf Proof:}

  In order to simplify the notation we call $E_n(f)= E_\mu(f\,|\,\mathcal{F}_n).$

  Suppose that $\psi$ is a KMS state.  Then for all
$a,b,c,d\in C(X)$ and all $n$ we have

\begin{equation*}
\psi((a e_n b)\sigma_{i\beta }(c e_n d)) =
  \psi((c e_n d)(a e_n b)).  \,\,\,\,\,\,\,\,(*3)
\end{equation*}

The left hand side  is equals to
  $$
  \psi(a \,e_n\, b \,\,c H^{-\beta[n]} e_n H^{\beta[n]} d) =  \psi(a S^n\, (S^*)^n\, b\,\, c H^{-\beta[n]} e_n H^{\beta[n]} d) =
  $$
  $$
  \psi(a E_n(b c H^{-\beta[n]}) \,e_n \,H^{\beta[n]} d) =
  \phi(a E_n(b c H^{-\beta[n]}) \,\lambda^{-[n]}\,  H^{\beta[n]} d).
  $$

  The right hand side of (*3) is equals to
  $$
  \psi(c E_n(d a) e_n b) =
  \phi(c E_n(d a) \lambda^{-[n]} b).
  $$

Now take $b=1$, $c = H^{\beta[n]}$, and $d = H^{-\beta[n]} \lambda^{[n]}$
and from (*3) we get
  $$
  \phi(a) =
  \phi( H^{\beta[n]} E_n( H^{-\beta[n]} \lambda^{[n]} a) \lambda^{-[n]}) =
  \phi(\Lambda^{-[n]} E_n(\Lambda^{[n]} a)).
  $$
  \medskip

  Now, we want to prove the other implication.

  Note that $\phi_\nu(ab)=\phi_\nu(ba)$ for continuous functions $a$ and $b$.

We would like to prove that
\begin{equation*}
  \psi((a e_n b)\sigma_{i\beta}(c e_m d)) =
  \psi((c e_m d)(a e_n b)), \,\,\,\,\,\,\,\,(*4)
\end{equation*}

for all $a,b,c,d\in A$ and $n,m\in \mathbb{N}$.

Suppose first the case $n\leq m$.

By the important property i) we get that the left hand side of (*4) is equals to
  $$
  \psi(a e_n b c H^{-\beta[m]} e_m H^{\beta[m]} d) =
  \psi(a E_n(b c H^{-\beta[m]}) e_m H^{\beta[m]} d) =
$$
$$  \phi(a E_n(b c H^{-\beta[m]}) \lambda^{-[m]} H^{\beta[m]} d) =
  \phi(E_n(b c H^{-\beta[m]}) H^{\beta[m]} \lambda^{-[m]} d a).
  $$

  Observe that
$H^{\beta[m]} (x)=  H^{\beta[n]} (x)\, H^{\beta[m-n]} (T^n(x))$ so the above is equals to
  $$
\phi(E_n(b c H^{-\beta[n]}H^{-\beta[m-n]} (T^n\,)  )\, H^{\beta[m-n]} (T^n)\, H^{\beta[n]} \lambda^{-[m]} d a) =
  $$
  $$
  \phi(E_n(b c H^{-\beta[n]}) H^{\beta[n]} \lambda^{-[m]} d a) =$$
  $$
  \phi(\Lambda^{-[n]} E_n(\Lambda^{[n]} E_n(b c H^{-\beta[n]})
H^{\beta[n]} \lambda^{-[m]} d a ) ) =$$
$$
  \phi(\Lambda^{-[n]} E_n(b c H^{-\beta[n]}) E_n(\lambda^{[n]}
\lambda^{-[m]} d a ) )=$$
$$ \phi(\Lambda^{-[n]} E_n(b c H^{\beta[n]})\, \lambda^{[n]}
\lambda^{-[m]}  E_n(d a ) ).
$$

  where in the last equality we use the fact that $\lambda^{-[m]} \lambda^{[n]} =\lambda^{-[m-n]}(T^n).$

By the other hand the right hand side of (*4) is  equals to
  $$
  \psi(c e_m E_n(d a) b) =
  \phi(c \lambda^{-[m]} E_n(d a) b) =
  \phi(b c \lambda^{-[m]} E_n(d a) ) =$$
  $$
  \phi(\Lambda^{-[n]} E_n(\Lambda^{[n]} b c \lambda^{-[m]} E_n(d a)
)) .$$
$$
  \phi(\Lambda^{-[n]} E_n(b c \lambda^{-[m]} \lambda^{[n]} H^{-\beta[n]} ) E_n(d a))=
  $$
  $$
  \phi(\Lambda^{-[n]} E_n(b c  H^{-\beta[n]} )\,\lambda^{-[m]} \lambda^{[n]}\, E_n(d a))
  $$

where in the last equality we use once more the fact that $\lambda^{-[m]} \lambda^{[n]} =\lambda^{-[m-n]}(T^n).$

In this way we showed the KMS condition in the case  $n\leq m$.
\medskip

For the case  $n\geq m$, using the important property j) we note that the left hand side of (*4)
is
  $$
  \psi(a e_n b\,\,\,\ c \,\,H^{-\beta[m]} e_m H^{\beta[m]} \,d) =$$
  $$
  \psi(a e_n E_m(b c H^{-\beta[m]}) H^{\beta[m]} d) =
  $$
  $$
  \phi(a \lambda^{-[n]}  E_m(b c H^{-\beta[m]}) H^{\beta[m]} d) =
  $$
 $$ \phi(
    \Lambda^{-[m]} E_m(\Lambda^{[m]}
    \lambda^{-[n]}  E_m(b c H^{-\beta[m]}) H^{\beta[m]} d a
    )) =$$
    $$
  \phi(\Lambda^{-[m]} E_m(b c H^{-\beta[m]}) E_m( H^{\beta[m]} d a
\Lambda^{[m]} \lambda^{-[n]} )) =$$
$$
  \phi(\Lambda^{-[m]} E_m(b c H^{-\beta[m]}) E_m( d a
\lambda^{[m]} \lambda^{-[n]} )).
  $$
  The right hand side of (*4) equals
  $$
  \psi(c E_m(d a)e_n b) =
  \phi(c E_m(d a)\lambda^{-[n]} b) =
  $$
  $$
  \phi(\lambda^{-[n]} b c E_m(d a) ) =
  \phi( \Lambda^{-[m]} E_m(\Lambda^{[m]} \lambda^{-[n]} b c E_m(d
a)) )=$$
$$
  \phi( \Lambda^{-[m]} E_m(b c H^{-\beta[m]} \lambda^{[m]} \lambda^{-[n]} ) E_m(d
a) ).
  $$
  The conclusion follows at once  because
$\lambda^{[m]} \lambda^{-[n]} \in \mathcal{F}_m$.

\qed

\bigskip

{\bf Corollary 2.1.} Suppose $\nu_\beta $ is an eigenprobability for the Ruelle operator of the potential
$ - \beta \log H.$ If the $C^*$-dynamical  system  state $\psi_{\nu_\beta} : \mathcal{U} \to \mathbb{C} $ is defined  by
$$ \psi_{\nu_\beta}(M_f e_n M_g) = \phi_{\nu_\beta} (G(M_f e_n M_g))\,=\, \int f \lambda^{-[n]} \,g\, d \nu_\beta ,  $$
then, $\psi_{\nu_\beta}$ is a KMS  state for temperature $\beta$.

{\bf Proof:} This follows from last theorem and Lemma 1.1

\qed

\medskip

Note that when $H$ is constant then $\mu$ is an eigenprobability for the associated Ruelle operator for any $\beta>0$.
From expression (*5) we can see that $\sigma_t$ in this case is the identity for any $t$. Moreover, by the KMS
relation $\psi_\mu (a\, b) = \psi_\mu(b\,a)$.

\medskip

We can ask about uniqueness of the KMS state.
To address this question  is the purpose of the next results.
\medskip

Our main theorem says:

{\bf Theorem 2.2:}
If $H$ is Holder positive and $\mu$ is a
Gibbs state for $p$ Holder, then for any given $\beta\in {\bf R}$,
a KMS  state $\psi$ in ${\cal U} (\mu)$
exists, it is unique and of the form
$$\psi_{\beta} (b) = \int \frac{f\, g}{\lambda^{[n]}} d \nu_{H,\beta},\, \, \, \forall b= M_f e^n M_g
\in {\cal U},$$
where $\nu_{\beta}$ is the eigenmeasure for
${\cal L}_{-\beta \log H}^{*}$.

{\bf Proof of Theorem 2.2:}

The existence of a  KMS follows  from the results from above. We fixed $\beta$.

Now we want
to show precisely how one can associate a Gibbs measure to a KMS state. We denote such KMS state by $\psi.$
We will denote $\psi_{\beta}$ the KMS state obtained from $\nu_\beta$.

Suppose $\psi$ is a KMS state, where the $H$ is fixed and  defines the semigroup $\sigma_t$.

Given the KMS state $\psi$, then $\psi (M_f)= \psi (f) $ defines a
continuous positive linear functional over $C(X)$ such
that $\psi (M_1)=1$. Therefore by Riesz Theorem,
there exists a probability $\nu$ such that  for any $f \in C(X)$
we have $\psi (f)= \int f d\nu= \int f S^0 (S^{*})^0 d \nu$.

The above definition takes in account just $n=0$ in a) above.
Remains the question: what conditions are imposed on
$\nu$ (defined from
$\psi$ as above) due to the fact that $\psi$ is a KMS state
for $H,\beta$?

This $\nu$ is our candidate to be the one associated to $\psi$ via  $\psi=\psi_\nu=\phi_\nu \circ G$ where hopefully $\nu$ satisfies
$$\int f\, d \nu = \int (\Lambda_n)^{-1} E_\mu( \Lambda_n\, f \,| \,\mathcal{F}_n\,) \,d \nu.
$$
for all continuous $f$, and also
$$ \psi_\nu(M_f e_n M_g) = \phi_\nu (G(M_f e_n M_g))\,=\, \int f \lambda^{-[n]} \,g\, d \nu .  $$

\medskip
Now we will show a recurrence relation which do not assume any  KMS state condition for $\psi$.
\medskip

{\bf Lemma 2.1:} Suppose that the $C^*$-state $\psi$ is such that $\psi(f)= \int f \, d\,\nu$, for any $f \in C(X)$.

Then, for any  $f \in C(X)$ and $n \in \mathbb{N}$
\begin{equation} \label{o1} \psi (f \,e^n)= \psi (f \,S^n (S^{*})^n)= \psi (S^n (S^{*})^n\,f)=\int f \lambda^{-[n]} d \nu\,.
\end{equation}

In other words, if $G (f \,S^n (S^{*})^n)= G (M_f e^n\, ) =  f \lambda^{-[n]}$ for any continuous function $f$, and $\phi_\nu(f)= \int f d \nu=\psi(f)$, then
$$\psi (f \,e^n) = \phi_\nu (G( f \, e^n)).$$

%Moreover, for any element $b\in \mathcal{U}$¨%\begin{equation} \label{o2}
%\psi (b \,S^n (S^{*})^n)=\,\psi( b\,\lambda^{-[n]}).
%\end{equation}

{\bf Proof:}

The first claim of the lemma will follow from
$$\psi (f S^n (S^{*})^n)= \psi (f\, (\lambda \circ T^n) \,  S^{n+1} (S^{*})^{n+1}) .$$

Indeed, for instance we get $\psi (f )= \psi (f\, \lambda  \,  S (S^{*})) $, and,
$\psi (f S\, (S^{*}))= \psi (f\, (\lambda \circ T) \,  S^{2} (S^{*})^{2}) $, and so on.

We need a preliminary estimate before proving the lemma.

For the transformation $T$, consider a partition $A_1,...,A_k$
of $X$ such that $T$ is injective in each $A_i$. Our
proof bellow is for the shift
in the Bernouilli space.
In the case
of the Bernoulli space with $k$ symbols $A_i$ is the cylinder
$\overline{i}$ with
first coordinate $i$. Now we consider
a partition of unity given by $k$ non-negative
functions $v_1,...,v_k$ such that
each $v_i(x)= I_{\overline i}(x)$ (the indicator function of the cilinder
$   {\overline i}$)  which has
support on $A_i$ and $\sum_{i=1}^k v_i (x) =1$ for all
$x\in X$. In the case $X$ is
the unitary circle and $T$ is expansive, using
a conjugacy with the shift,
we obtain similar results.

Denote now the functions $u_i$
given by: if $x $ is in the cylinder $\overline{i}$
\begin{equation*}  u_i(x) = ( \lambda(x))^{1/2}
=( p(x)^{-1})^{1/2},\,\,\,\,\,\,\,\,(*1)
\end{equation*}
for each $i\in \{1,...,k\}$,
so, $\sum_{i=1}^k u_i^2 (x) =\lambda(x)= p^{-1}(x)$, for any $x\in X$.

An easy computation shows that $\sum_{i=1}^k M_{u_i} S S^{*}
M_{u_i} = 1$. More precisely, if $ x$ is in the cylinder $\overline{i}$,
then given $\eta$, we have
\begin{equation*} [\,M_{u_i} S S^{*} M_{u_i} (\eta)\,] (x) =   u_i(x) \, \sum_{\{z\, |\,\sigma(z) =
\sigma (x)\}}
 p(z) \, u_i(z) \, \eta(z)\,=\,\eta(x).\,\,\,\,\,\,\,\, (*2)
 \end{equation*}

 Indeed,
for $x$ in the cylinder $\overline{i},$
take  $u_i (x)= p^{-1/2} (x).$ This is so because for $x=(j,x_2,x_3,...)$ we get
$$\sum_{i=1}^k M_{u_i} S S^{*} M_{u_i} (\xi)(x)=  M_{u_j} S S^{*} M_{u_j}(x) =p^{-1/2} (x) [\mathcal{L}_{p}  (\xi\, u_j) ] (\sigma(x))=$$
$$p^{-1/2} (x)\sum_{i=1}^d  p(i,x_2,x_3,..)  u_j (i,x_2,x_3,..)  \xi (i,x_2,x_3,..)=$$
$$p^{-1/2} (x)  p(j,x_2,x_3,..)  u_j (j,x_2,x_3,..)  \xi (j,x_2,x_3,..)=$$
$$p^{-1/2} (x)  p(x) p^{-1/2} (x)  \xi (x)=\xi(x).$$

Now, we continue the argument:

$$S^n (S^{*})^n= S^n  1 (S^{*})^n=
S^n [\,\sum_{i=1}^k   (M_{u_i} S S^{*} M_{u_i}) \,] (S^{*})^n=$$
$$
\sum_{i=1}^k (S^n  \, (M_{u_i} S S^{*} M_{u_i}) \,  (S^{*})^n).$$

Now we use the relations $S^n M_g= M_{\alpha^n (g)} S^n$
and $M_g (S^{*})^n = (S^{*})^n M_{\alpha^n (g) }$ in last
expression and we get

$$S^n (S^{*})^n=
\sum_{i=1}^k
\, M_{\alpha^n (u_i)} S^n \, S  S^{*} \,  (S^{*})^n M_{\alpha^n
(u_i)} \, =
\sum_{i=1}^k
\, M_{\alpha^n (u_i)} S^{n+1}  (S^{*})^{n+1} M_{\alpha^n
(u_i)} \, $$

Now we will prove the lemma.

Using last expression and then Remark 2 for
$g=\alpha^n (u_i)\in C(X)$ and $a= S^{n+1}  (S^{*})^{n+1}$  we get

$$ \psi ( M_f S^{n}  (S^{*})^{n}) = \psi (M_f \, \sum_{i=1}^k
M_{\alpha^n (u_i)} S^{n+1}  (S^{*})^{n+1} M_{\alpha^n
(u_i)} \, )=$$

$$  \psi (M_f \, \sum_{i=1}^k
M_{\alpha^n (u_i)}M_{\alpha^n
(u_i)} S^{n+1}  (S^{*})^{n+1}  \,)
=$$
$$\psi (M_f \, \sum_{i=1}^k
M_{\alpha^n (u_i)^2}  S^{n+1}  (S^{*})^{n+1}  \,  )
=$$
$$\psi (M_f \,
M_{\alpha^n (\sum_{i=1}^k  (u_i)^2)}  S^{n+1}  (S^{*})^{n+1}  \,  ).
=$$
$$\psi (M_f \,  (\lambda\circ T^n) S^{n+1}  (S^{*})^{n+1}   )
$$

This shows the claim of the lemma.

%Note that the same proof as above shows that for any $b$ in the $C^*$-Algebra $\mathcal{U}$ we get

%$$  \psi (b \, S^{n}  (S^{*})^{n}  \,)
%= \psi (b \,  (\lambda\circ T^n) S^{n+1}  (S^{*})^{n+1}   ) .$$

%From this follows the second part of the lemma.

\qed

We denote $E_m (f) = E_\mu(f\,|\, \mathcal{F}_m)$.

{\bf Corollary 2.1} If $\psi$ is Gibbs for $H$ at temperature zero, and $\nu$ is such that for any continuous function $f$ we have $\psi(f)= \int f d \nu$, then

$$\phi_\nu  (f\,) = \phi_\nu (\,\,(\Lambda_n)^{-1} E_\mu( \Lambda_n\, f \,| \,\mathcal{F}_n\,) \,), $$

which is the same that to say that $\nu$ satisfies

\begin{equation*} \int f\, d \nu = \int (\Lambda_n)^{-1} E_\mu( \Lambda_n\, f \,| \,\mathcal{F}_n\,) \,d \nu.  \,\,\,\,\,\,\,\,(*6)\end{equation*}

{\bf Proof:} We get from
last lemma that $\psi (f \,e^n) = \phi_\nu (G( f \, e^n))$ where
$\phi_\nu(f)= \int f d \nu=\psi(f)$. Now, from Theorem 2.1 we get that (*6) is true.

\medskip

\qed

\medskip
Now we will show the uniqueness of the KMS state:

\medskip

{\bf Theorem 2.3:} Given any KMS $\psi$,
then $\psi=\psi_\beta$ where $\psi_\beta$ is the KMS state associated to
the Gibbs probability $\nu_\beta$.

{\bf Proof:}

In order to do that we will show that any possible $\nu$ as defined above from the KMS $\psi$
is equal to $\nu_\beta$.

Take $\nu$ a probability associated to $ \psi$, then
for each $n$, and $f\in C(X)$ we have
\begin{equation*} \int f d\nu = \int E_n (f \Lambda_n^{-1}) \Lambda_n d \nu= \int \alpha^n
({\cal L}^n_{\beta} (f ) ) \Lambda_n d \nu. \, \, \,\,\,\,\,\,\,\,\,(*7)
\end{equation*}

We claim that
$$ \lim_{n\to \infty}  \int E_n (f \Lambda_n^{-1}) \Lambda_n d \nu=
\int f d \nu_\beta,$$
and this shows that $\nu=\nu_{\beta}$, and therefore
$\psi=\psi_\beta$.

Now we show the claim. Note that
$$\int f d \nu=    \int E_n (f \Lambda_n^{-1}) \Lambda_n d \nu=
\int \alpha^n ({\cal L}^n_{\beta} (f )) \Lambda_n d \nu =
\int \alpha^n (\frac{{\cal L}^n_{\beta} (f )}{\lambda_\beta^n}) \Lambda_n
\lambda_\beta^n d \nu,$$
where $\lambda_\beta$ is the eigenvalue associated to ${\cal L}_{\beta} $.

Applying the above expression to $f=h_\beta$ (we can assume $h_\beta$ is such that
$\int h_\beta d \nu_\beta=1$)
and
using the fact that  $ {\cal L}^n_{\beta}
(h_\beta )= \lambda_\beta^n h_\beta $ we get
$$ 0<d = \int h_\beta d\nu = \int \alpha^n  ( h_\beta)  \Lambda_n \lambda_\beta^n  d \nu.$$

As $h_\beta$ is continuous and positive, there exists $c>0$ such
for all $x\in X$ we have $h_\beta(x)>c$.

From this follow that
$$d= \int   \alpha^n  (h_\beta)  \Lambda_n \lambda_\beta^n d \nu >
c \int  \lambda_\beta^n
 \Lambda_n d \nu.$$

Therefore,
$$  \int  \lambda_\beta^n
\Lambda_n d \nu < d/c$$

Denote $I=\int f d \nu_\beta.$

It is known (see [Bo]) that uniformly in  $z\in X$,
we have
$$\lim_{n \to \infty}
\frac{  {\cal L}^n_{\beta} (f )(z)}{\lambda_\beta^n}= h_{\beta} (z) I=
h_{\beta} (z) \int f d \nu_\beta .$$

Therefore, given $\epsilon>0$, we can find $N>0$
such that for all $n>N$ we
have for all $z\in X$
$$| \frac{  {\cal L}^n_{\beta} (f )(z)}{\lambda_\beta^n}- I h_\beta (z)|
\le \epsilon.$$

Then, for $n>N$
$$ |\int
\frac{ \alpha^n ({\cal L}^n_{\beta} (f ))}{\lambda_\beta^n} \Lambda_n
\Lambda_\beta^n d \nu - \int I \alpha^n (h_\beta )  \Lambda_n
\Lambda_\beta^n d \nu \, |\leq $$
$$\int | \frac{ \alpha^n ({\cal L}^n_{\beta} (f ))}{\lambda_\beta^n} (y) -
I \alpha^n (h_\beta)  (y)| \Lambda_n (y)
\Lambda_\beta^n (y) d \nu=$$
$$ \int | \frac{  {\cal L}^n_{\beta} (f )}{\lambda_\beta^n}
(T^n (y)) -
I h_\beta  (T^n (y))| \Lambda_n (y)
\lambda_\beta^n (y) d \nu \leq \frac{\epsilon d}{c}$$

The conclusion from (*7) is that for any $f\in C(X)$
$$   \lim_{n\to \infty}  I\, \int  \alpha^n (h_\beta ) \Lambda_n \lambda_\beta^n d \nu =
\int f d \nu  .$$

Consider now $f=1$ and we get
$$  \lim_{n\to \infty}  \int  \alpha^n (h_\beta )
\Lambda_n \lambda_\beta^n d \nu = 1.$$

From this we conclude that $\int f d\nu=I= \int f d \nu_\beta$ for all
$f\in C(X)$.

This shows the uniqueness and that $\nu=\nu_\beta.$

$\square$

The final conclusion is that any KMS $\psi$ for $H,\beta$ is equal to the
$\psi_\beta$ associated to $\nu_\beta$.

\vspace{1.0cm}

\centerline {\bf Section 3 - no phase transitions} \vspace{1.0cm}

We consider here an interesting example of a KMS state associated with the reference measure $\mu$
given by the maximal entropy measure for the shift
in $2$ symbols $\{0,1\}$. In this case $p=1/2$ is contant.
We will define a special  potential $H$ and we will consider
specifically   the special value $\beta=1$

We refer the reader to [H] [L2] [L3] [FL] [Y] [L]  for references and results about the topics
discussed in this section.

We are going to introduce the Fisher-Fedenhorf model of Statistical Mechanics
in the therminology of Bernoulii spaces and Thermodynamic Formalism [H].

We define  $\Sigma^+$ to be the shift space $\Sigma^+= \Pi_0^\infty\{ 0,1\}$
and denote by $T: \ \Sigma^+ \to \Sigma^+$  the left shift map.
We write $z=(z_0 z_1\dots)$ for a point in $\Sigma^+$
and $[w_0 w_1 \dots w_k]= \{ z:\ z_0=w_0, z_1=w_1, \dots z_k=w_k\}$ for a
cylinder set of $\Sigma^+$.

%We write $I_A$ for the indicator function of a set $A$; thus
%$I_A(z)=1$ if $z\in A$ and $=0$ if $z \notin A$.

We denote by $M_k \subset \Sigma^+, $ for $ k>1,$ the
cylinder set
$[\underbrace{111\dots 11}_{k} 0]$ and by $M_0$  the cylinder set
$[0].$  The ordered collection  $(M_k)_{k=0}^\infty$
is a partition of $\Sigma^+$; in other words these
 sets are disjoint and their
union is the whole space (minus the point $(11\dots)$).
Note that $T$ maps $M_k$ bijectively onto $M_{k-1}$ for $k\geq 1$,
and onto $\Sigma^+$ for $k=0$.

The point $(1111...)$ is fixed for $T$.

For $\gamma>1$  a fixed real constant,
we consider the   potential $g(x)$ such that $g(111111\dots .)=0$,
$$
g(x)=a_k=-\gamma\log\biggl(\frac{k+1}{k}\biggr),$$
for
$x\in M_k$, for $k\neq 0,$  and
$$
a_0=-\log(\zeta(\gamma)),$$
for $x\in M_0,$
where $\zeta$ is the Riemann zeta function.

By definition,
$$
\zeta(\gamma)=(1^{-\gamma}+2^{-\gamma}+\dots)$$
and so the reason for defining $a_0$ in such way is that,
if we define  $s_k=a_0 + a_1 +\dots  + a_k$, then
$\Sigma e^{s_k}=1.$

From now on we assume $\gamma>2$, otherwise we have to consider
sigma-finite measures  and not probabilities in our problem.

The potential $1<(\frac{k+1}{k})^\gamma = H(x)= e^{-g(x)},$ for $x \in M_k,$
is not H\"older and in fact is not of summable variation. Note that
$H( 1111...)=1$, The pressure  $P(-\log H) =P(g)=P( \log p + \log 2- 1\,  \log H)=0$ and one can show that there exist two equilibrium states for such a potential $g$ (in the sense of  minimizing measures
for the variational problem):
a point mass (the Dirac delta $\delta{(111...)}$)  at $(1111\dots )$, and a second measure
which we shall denote by $\tilde{\mu}$ (see [H])

The existence of two probabilities $\tilde{\mu}$ and
$\delta_{(1111...)}$ for the variational problem of pressure
defines what is called a phase transition
in the sense of Statistical Mechanics [H] [L3].

We will describe bellow how to define this measure $\tilde{\mu}$.

Consider as in [H] ${\cal  L}_g^*$, the dual of the
Ruelle-Perron-Frobenius operator $ {\cal  L}_g$ associated to $g$,
where the action of $ {\cal  L}_g$
 on continuous functions is given by
$$  {\cal  L}_{\beta=1} (\psi)(y)=\sum_{T(x)=y}e^{g(x)}\psi(x).$$

The function $P(-\beta \log H) = P(\beta  g)$ is strictly
monotone for $\beta<1$ and constant equal zero for $\beta>1$ [H].

We claim that there is a unique probability measure $\nu$ on
$\Sigma^+$ which satisfies $ {\cal  L}_g^*\nu=\nu$ [FL] [H]. To
prove this, note first that $\nu$ cannot have any mass at
$(11\dots)$; it follows that $M_0$ has positive mass, and the
stipulation that $\nu$ be an eigenmeasure then gives a recurrence
relation  for the masses of $M_k$. Since $T(M_k)=M_{k-1}$ for
$k\geq 1$, we have that the masses of the sets in this partition
are
 $$\nu(k)= \nu(M_k)=e^{s_k}=\frac{(k+1)^{-
\gamma}}{\zeta(\gamma)},\, k \geq 0 ;$$
in particular,
$$\nu(0)=\nu(M_0)=e^{s_0}=e^{a_0}=\frac{1}{\zeta(\gamma)}.$$
By the same reasoning, $\nu$ is determined on all
higher cylinder sets for the partition  $(M_k)_{k=0}^\infty$.
Hence  $\nu $ exists and is unique.

The measure $\nu$ defined above is the unique eigenmeasure for
${\cal  L}_{\beta=1}^*$ and denoted by $\nu_1$.

The measure defined by the delta-Dirac on $(111...)$ is invariant
but is not a fixed eigenmeasure for ${\cal  L}_g^*$.

This measure $\nu_1$ defines a KMS state
$\psi_{\nu_1}$ for such $H$, $\beta=1$ and  ${\cal U} (\mu)$.

We conjecture that there is another KMS state $\psi$ different from $\psi_{\nu_1}$
but not associated to a measure.
Note that such $H$ assumes the value $1$ in just one point.

We   define $\tilde{h}(x)$ for $x\in M_t$  by
$$\tilde{h}_t= \tilde{h}(x)=\nu(t)^{-
1}\sum^\infty_{i=t}\nu(i).$$

The function $\tilde{h}$ satisfies $ {\cal  L}_g
(\tilde{h})=\tilde{h}.$

The integral $\int \tilde{h} (x) d \nu_{1} (x)$ is finite if and only if
$\gamma>2$. One can normalize $\tilde{h}$,
multiplying by a constant $u$ to get
 $h= u  \tilde{h}$ with $\int h d \nu_{1} =1.$

This constant is
$$u= \frac{1}{ \sum_{t=1}^{\infty}  t \nu(t-1) }=  \frac{\zeta (\gamma)}{ \sum_{t=1}^{\infty}  t^{1-\gamma} }= \frac{\zeta(\gamma)}{\zeta (\gamma - 1)}.$$

 The probability $\tilde{\mu}$ has positive entropy  and
its support is all $\Sigma^{+}$ (see [H] or [L3] [FL]).

Consider now the invariant probability measure $\tilde{\mu}=h\nu_{1}.$ It is known
that  $\tilde{\mu}$ is
an equilibrium state for $-\log H$ in the variational sense ($\beta=1$) [H].
It is easy to see (because $-\log H(11111..)=-\log 1=0$) that the Dirac-delta measure $\delta_{(11111...)}$
is also an equilibrium state for $-\log H$ in the variational sense ($\beta=1$).

The probability $\tilde{\mu}$ has positive entropy  and
its support is all $\Sigma^{+}$ (see [H] or [L3] [FL]).

We can conclude from the above considerations that not always an
equilibrium probability $\rho$ for the
pressure is associated to a KMS state $\psi_\rho$ whitout the hypothesis
of $H$ and $p$ been Holder. In the present example, this happen
because $\rho=\delta_{(1111...)}$ is not
an eigenmeasure of the dual of the Ruelle-Perron-Frobenius operator ${\cal L}_\beta$
but it is an equlibrium measure for $\beta=1$.

In [L2] and [L3] the lack of differentiability of the Free energy
is analyzed  and in  [L3] [Fl] [Y] it is shown that such systems
present polynomial decay of correlation. In [L1]
it is presented a dynamical model with three equilibrium states.

\vspace{1.0cm}

\centerline {\bf Bibliography}
\vspace{1.0cm}

[Bi] P. Billingsley
{\it Probability and Measure}, Wiley, New York, (1968)

[Bo] R. Bowen, Equilibrium states and the ergodic theory
of
Anosov diffeomorphisms, {\it Lecture Notes in Mathematics} 470,
Springer Verlag, (1975)

[BS] R. Bowen and C. Series, Markov maps associated with a Fuchsian group,
{\it IHES Publ Math}, N 50, pp 153-170 (1979)

[BR] O. Bratelli and W. Robinson, {\it
Operator Algebras and Quantum Statistical Mechanics},
Springer Verlag, (1994)

[GL] G. Castro and A. Lopes, {\it KMS States, Entropy and a Variational Principle for Pressure}, Real Analysis Exchange, Vol 34 , Issue 2, 333-346 (2009)

[CL1] L. Cioletti and A. O.  Lopes,
Interactions, Specifications, DLR probabilities and the
Ruelle Operator in the One-Dimensional Lattice,
Discrete and Cont. Dyn. Syst. - Series A, Vol 37, Number 12,
6139-6152 (2017)

[E1] R. Exel, Crossed-Products by finite Index Endomorphisms and
KMS states,  {\it J. Funct. Anal.} 199, no. 1, pp 153-188 (2003).

[E2] R. Exel, A new look at the Crossed-Products of a
$C^{*}$-Algebra by endomorphism, {\it  Erg. Theo. Dyn. Syst.} 23,
no. 6, pp 1733--1750 (2003)

[E3] R. Exel, KMS states for generalized Gauge actions on
Cuntz-Krieger Algebras, {\it Bull. Braz. Math. Soc.} 35, n. 1, pp
1-12,  (2004)

[EL1] R. Exel and M. Lacca, Cuntz-Krieger Algebras for infinite matrices,
{\it J. reine angew. Math}, 512, pp 119-172, (1999)

[EL2] R. Exel and M. Lacca, Partial Dynamical Systems and the KMS
condition, {\it Comm. Math. Phys.} 232, no. 2, pp 223-277 (2003).

[EL3] R. Exel and A. Lopes, $C^*$ Algebras, approximately proper
equivalence relations and Thermodynamic Formalism, {\it Erg. Theo.
and Dyn. Syst.}, Vol 24, pp 1051-1082 (2004)

[E0] R. Exel,
Inverse semigroups and combinatorial $C^*$-Algebras, Arxiv (2008)

[FL] A. Fisher and A. Lopes, Exact bounds for the polynomial decay of correlation,
1/f noise and the central
limit theorem for a non-Holder Potential, {\it Nonlinearity}, 14, pp 1071-1104 (2001).

[FF] B. Felderhof and M. Fisher, four articles, {\it Ann. Phys. 58}, pp 176-281 (1970)

[H] F. Hofbauer, Examples for the non-uniqueness of the equilibrium
states, {\it Trans. AMS}, 228, pp 133-149 (1977)

[J] V. Jones, Index for subfactors, {\it Invent Math}, 72, pp 1-25, (1983)

[K] H. Kosaki, Extensions of Jones' Theory on index to arbitrary subfactors,
{\it J. Funct Anal}, 66, pp 123-140, (1986)

[L] M. Lacca, Semigroups of *-endomorphisms, Dirichlet Series and phase transitions,
{\it J. Funct Anal}, 152, pp 330-378, (1998)

[L1] A. Lopes, Dynamics of Real Polynomials on the
Plane and Triple Point Phase Transition,
{\it Mathematical and Computer Modelling\/}, Vol. 13,
N$^{\underline o}$ 9, pp. 17-32, 1990.

[L2] A. Lopes, A First-Order Level-2 Phase Transition
in Thermodynamic Formalism, {\it Journal of Statistical
Physics\/}, Vol. 60, N$^{\underline o}$ 3/4, pp.
395-411, 1990.

[L3] A. Lopes, The Zeta Function, Non-Differentiability
of Pressure and The Critical Exponent of Transition",
{\it Advances in Mathematics\/}, Vol. 101, pp. 133-167,
1993.

[L4] A. Lopes, Entropy, Pressure and Large Deviation, {\it Cellular Automata, Dynamical Systems and Neural Networks\/},
E. Goles e S. Martinez (eds.), Kluwer, Massachusets, pp. 79-146, 1994.

\end{document}